\documentclass[11pt]{amsart}

\usepackage{amssymb}
\usepackage{amscd}
\usepackage{a4}
\usepackage{amsthm}

\usepackage{fullpage}

\usepackage{etex}
\usepackage{dsfont}
\usepackage[matrix, arrow, curve, frame, arc]{xy}
\usepackage{pgf}
\usepackage{tikz}
\usetikzlibrary{arrows,shapes,automata,backgrounds,petri,calc}

\newcommand{\tikzAngleOfLine}{\tikz@AngleOfLine}
\def\tikz@AngleOfLine(#1)(#2)#3{%
\pgfmathanglebetweenpoints{%
\pgfpointanchor{#1}{center}}{%
\pgfpointanchor{#2}{center}}
\pgfmathsetmacro{#3}{\pgfmathresult}%
}

\ifx\textheight\undefined
\else
\advance\voffset by0.8cm
\advance\voffset by-2.3cm
\fi

\newtheorem{theorem}{Theorem}[section] 
\newtheorem{lemma}[theorem]{Lemma}

\newtheorem{main-theorem}[theorem]{Theorem}

\newtheorem*{problem*}{Problem}

\theoremstyle{definition}

\newtheorem*{question*}{Question}

\newcommand{\La}{\Lambda}

\newcommand{\bZ}{\mathbb{Z}}

\newcommand{\ul}{\underline}

\begin{document}

\title{A note on representation-finite symmetric algebras}

\author{Karin Erdmann}

\address{Mathematical Institute, University of Oxford,
ROQ, OX2 6GG, Oxford UK}

\email{erdmann@maths.ox.ac.uk}

\begin{abstract} We give  characterisations for representation-finite 
symmetric algebras of period four, and describe their basic algebras.
In particular if such an algebra is indecomposable, it has at most
two simple modules.
\end{abstract}

\maketitle

\section{Introduction} 

Tame symmetric algebras of period four occur as 2-blocks of group algebras with 
quaternion defect groups, and algebras of quaternion type \cite{E},  and in much more generality,  as weighted surface algebras and their
virtual mutations and some other related algebras ( \cite{ES1} (and \cite{ES-GV}, \cite{HSS}, and \cite{SS}). 
In  \cite{ES3}  the main result is the classification of  tame symmetric algebras of period four whose Gabriel quiver is 2-regular, and 
there is work in progress towards a general classification of tame symmetric algebras of period four.
As part of this project, we wish to clarify 
when exactly a symmetric algebra of period four is of finite type.
We  determine these algebras explioictly, there are three families and for each, the number of simples is at most three.  On the way we  give some characterisations, which are used in \cite{EHS}..

Recall that the syzygy $\Omega(M)$ of a module $M$ is the kernel
of a minimal projective cover of $M$. The module $M$ is periodic if 
$\Omega^n(M)\cong M$ for some $n\geq 1$, and if so, the minimal such $n$
is the period of $M$.
For general background we refer to \cite{ASS}.

\bigskip
\section{The characterisations}

We fix an algebraically closed field $K$. Since the properties we will study are invariant under Morita equivalence,
we can work with basic algebras, of the form $KQ/I$ where $Q$ is a finite quiver which we can take to be connected, and $I$ is an admissible ideal in 
the path algebra $KQ$. 

For a vertex $i$ of $Q$, we write $i^-$ for the set of arrows ending at $i$,
and $i^+$ for the set of arrows starting at $i$. We work with right
modules. We denote by $P_i$ the indecomposable projective module corresponding
to the vertex $i$, and by $S_i$ its simple quotient. Let also $p_i:= \ul{\dim} P_i$ its dimension vector. 

We assume each simple module $S_i$ has
$\Omega$-period diving $4$. Then 
there is an exact sequence
$$0 \to S_i \to P_i \stackrel{d_3}\to P_i^- \stackrel{d_2}\to P_i^+ \stackrel{d_1}\to P_i \to S_i \to 0
$$
with ${\rm Im} (d_k) \cong \Omega^k(S_i)$ for $1\leq k\leq 3$. 
Here $P_i^+ = \oplus_{\alpha: i\to j} P_j$ and $P_i^- = \oplus_{\gamma: k\to i} P_k$.

\medskip

\begin{lemma} \label{lem:2.1}
Let $\La$ be an indecomposable basic symmetric algebra such that all indecomposable non-projective modules have period dividing $4$. Then the following are equivalent:\\
(a)  There is a  simple module $S$ such that $\Omega^2(S)$ is simple.\\
(b) There is a vertex $i$ of the quiver with $|P_i| = |P_i^{+}|  = 
	|P_i^{-}|$ where $|-|$ denotes dimension.\\
(c) The Gabriel quiver of $\La$ is either one vertex with one loop, or it is of the form
 \[
 \xymatrix{
  1 \ar@<+.5ex>[r]^{\alpha}
   & 2 \ar@<+.5ex>[l]^{\beta} 
 }
\].\\
(d) The quiver of $\La$ has an arrow $\alpha: i\to j$ such that $i^+ = \{ \alpha \} = j^-$.\\
If  these hold then  $\La$ has finite type.
\end{lemma}

{\it Proof } 
(a) $\Rightarrow$ (b) \  Assume $S_i$ and $S_j$ are simple and $\Omega^2(S_i)\cong S_j$ (possibly $S_i\cong S_j$).  Then the right part of the exact sequence gives  
$0\to S_j\to P_i^+\to P_i\to S_i\to 0$ which implies $P_i^+|=|P_i|$. Similarly the left half implies $|P_i^-|=|P_i|$.

(b) $\Rightarrow$ (c) \ It follows from (b) that $\dim \Omega^2(S_i)=1$ and $S_j:=\Omega^2(S_i)$ is simple. Then $P_i^+=P_j$ and thee is a unique arrow  starting at $i$, and it ends at $j$. Denote it by $\alpha$. As well
$P_i^-\cong P_j$, so there is a unique arrow, $\beta$ say, ending at $i$ and it starts at $j$. 

Now $\Omega^2(S_j)\cong S_i$ since $\Omega^4(M)\cong M$ for all non-projective $M$. So we may interchange the roles of $i, j$. So there is a unique arrow starting at $j$, and there is a unique arrow ending at $j$. These are then $\beta$ and $\alpha$. 
Since $Q$ is connected, part (c) follows.
The implication (c) $\Rightarrow$ (d) is obvious.

(d) $\Rightarrow$ (a)  If (d) holds then $\Omega(S_i)= \alpha\La$ and $\Omega^{-1}(S_j)\cong \alpha\La$. It follows that $\Omega^2(S_i)\cong S_j$.

Any of the conditions above  implies $\La$ has finite type:  Part (c) shows that 
 $\La$ is a Nakayama algebra, which is   finite type (see for example
\cite{SY} Vol I, Ch I.10.3). Alternatively, it is easy to determine a presentation of a symmetric algebra with quiver as in (c).
$\Box$

\bigskip

\section{The basic algebras}

We describe now  which symmetric algebras
of finite type have only indecomposable 
modules of $\Omega$-period dividing four.
The algebras appearing in the previous Lemma appear in (i) and (ii),  and there
is only one other family of algebras.

\bigskip

	\begin{lemma} \label{lem:3.1} Assume $\La$ is symmetric indecomposable of finite type such that
 all indecomposable non-projective
modules have $\Omega$-period dividing $4$. Then $\La$ is one of the following.\\
(i) \ $\La \cong K[T]/(T^n)$ for some $n\geq 1$. \\
(ii) \ $\La \cong KQ/I$ where $Q$ is the quiver 
\[
 \xymatrix{
  1 \ar@<+.5ex>[r]^{\alpha}
   & 2 \ar@<+.5ex>[l]^{\beta} 
 }
\]
and $I = \langle (\alpha\beta)^n\alpha, (\beta\alpha)^n\beta \rangle$ for some $n\geq 1$.\\
(iii) \  $\La = KQ/I$ where $Q$ is the quiver 
\[
 \xymatrix{
  1 \ar@(dl,ul)[]^{\rho} \ar@<+.5ex>[r]^{\alpha}
	& 2  \ar@<+.5ex>[l]^{\beta} 
}
\]
and 
	$I = \langle  \rho\alpha, \beta\rho,  (\alpha\beta) -\rho^t\rangle$ for some 
$t\geq 2$.
\end{lemma}

\medskip

We will exploit the stable AR quiver.
For a symmetric algebra of finite type, by \cite{R},
the stable AR quiver isomorphic to
$\bZ\Delta/G$ where $\Delta$ is
a Dynkin diagram of type ADE, and $G$ is an admissible group of automorphisms
of $\bZ\Delta$.
Our main ingredient is a result due to Butler and Ringel \cite{BR}:

	\begin{lemma}\label{lem:3.2} Assume $\La$ is a basic algebra, and $\alpha: i\to j$ is an arrow in the quiver $Q$ of $\La$. Then the AR sequence of the module $e_i\La/\alpha\La$ has indecomposable middle term.
	\end{lemma}

We say that $M$ is a Butler-Ringel module if it is indecomposable and
its AR sequence has indecomposable middle term.

When $\Delta$ has type $A_n$ then there are precisely two $\tau$-orbits
with Butler-Ringel modules. For types $D$ and $E$ there
are three such orbits.

The algebra $\La$ is symmetric, hence
the action of $\Omega$ induces an equivalence of the stable category, taking
AR sequences to AR sequences. Hence the action  induces
a permutation, $\omega$ say,  on the orbit graph. (Recall the orbit graph  has vertices the $\tau$-orbits,
and there is an edge $a - b$ precisely if there are $x, y$ in the
AR quiver such that $x\in a$ and $y\in b$ and either $x\to y$ or $y\to x$
in the stable AR quiver.) Since $\Omega^2 \cong \tau$ we have $\omega^2 = 1$.

\bigskip

{\it Proof of Lemma \ref{lem:3.1}} \  If $\La$ is a local algebra then $\La$ is as in part (i). Furthermore,
if $\La$ is not local but satisfies (b) of Lemma \ref{lem:2.1} then it is
the algebra in part (ii). So we assume now that $Q$ does not have an arrow as in(b) of Lemma \ref{lem:2.1}. 
In particular any  simple module of $\La$ must have $\Omega$-period four.
	We must show that $\La$ is an algebra as in (iii).

\bigskip

{\it (1) \ We claim that the Gabriel quiver $Q$ must have a vertex at which only
one arrow starts, or only one arrow ends}

{\it Proof } 
To have finite type, the separated quiver of $\La$ must be a disjoint union of
Dynkin quivers. This means that there are at most two  arrows starting
or ending at a given vertex. Then there must be a vertex at which only one
arrow starts, or ends (otherwise the separated quiver would have components
of the form $\widetilde{A}_n$). 

\bigskip

Say $\alpha: 1\to 2$ and $\alpha$ is the only arrow ending at vertex $2$.
There must be an arrow starting at $2$, $\beta$ say (the quiver of a symmetric algebra does not
have sinks).
Furthermore, 
since we exclude Lemma 2.1(b) there must be
an arrow $\rho \neq \alpha$ starting at $1$. 

Consider the $\Omega$-orbit of the simple module $S_2$. 
We have that  $\Omega^{-1}(S_2)$ is isomorphic to $\alpha\La$. This is therefore
a Butler-Ringel module, 
 and therefore all syzygies are Butler-Ringel modules.
This gives four indecomposables,
$$S_2, \Omega(S_2), \Omega^2(S_2), \ \alpha\La 
\leqno{(*)}$$
(2) {\it We claim that the stable AR quiver  cannot be of type $D$ or $E$:}
Suppose we had type $D$. Then the two $\tau$-orbits of length two where the
AR sequences have indecomposable middle terms, must correspond to the branch 
vertices in the orbit graph. But then $\Omega$ fixes all other $\tau$-orbits, 
that is modules other than the modules in (*) have $\Omega$-period $2$. However,
by our assumption, the simple module $S_1$ has $\Omega$-period four, a contradiction.

Suppose we had type $E$, we may label the vertices of the orbit graph in order
as $1, 2, 3, 5, 6$ where $3$ is the branch vertex, and the last vertex is $4$. Then 
 $\omega = (1 6)(2 \ 5)(3)(4)$ and 
the  modules in (*) form the $\tau$-orbits corresponding to $1, 6$. The modules corresponding to 3 and 4 are fixed by $\omega$ and it follows that the modules in these two $\tau$-orbits have $\Omega$-period $2$. 
Consider the position of the simple module $S_1$. It is one of the four modules
corresponding to $2$ and $5$. Let $X$ be one of these
 four modules, then  there is either an irreducible map $S_2\to X$ or $X\to S_2$, or there is either an irreducible map $X\to \Omega(S_2)$ or $\Omega^{-1}S_2 \to X$. In all cases, $X$ cannot be isomorphic to $S_1$. 
It follows that $S_1$ must be one of the modules corresponding to $3$ or $4$ and then $\Omega^2(S_1)\cong S_1$, a contradiction.

\bigskip

Hence the stable AR-quiver is type $A$, and the four modules (*) are precisely 
the Butler-Ringel modules.
Therefore  
the arrow modules $\beta\La$ and $\rho\La$ must occur in (*), by Lemma 
\ref{lem:3.2}.  It follows that  
$\beta$ ends at $1$, and $\beta\La = \Omega(S_2)$, and moreover 
$\rho\La \cong \Omega^2(S_2)$ and $\rho$ is a loop at vertex $1$.
This implies that $\rho\La = \Omega(\beta\La)$ and $\Omega(\rho\La) = \alpha\La$. In other words, we have
$$\beta\rho=0, \ \ \rho\alpha = 0.
$$
Suppose $\alpha\beta$ does not occur in some minimal relation. It is nonzero in $\La$ since $\La$ is symmetric.
Consider the  idempotent algebra $e_1\La e_1$. Then $\rho$ is an arrow for this
algebra. If $\alpha\beta$ does not occur in some minimal relation then
it is another loop in $e_1\La e_1$ and $e_1\La e_1$ is not of finite type.
It follows that 
$\La$ is not of finite type, a contradiction.

Then we have a minimal relation $\alpha\beta = \rho^t\cdot u$ where $u$ is a unit in $e_1\La e_1$, that is it is a polynomial in $\rho$ with nonzero constant term and we may assume the constant term is $1$.  
The given information implies that the image of 
$\alpha\beta$ is in the socle of $\La$.  
Then $\alpha\beta = \alpha\beta u^{-1} = \rho^t$ 
and we have 
the stated presentation.

\end{document}